\documentclass[a4paper,12pt]{article}

\usepackage[english]{babel}
\usepackage{amsfonts,amsmath,amsthm,mathrsfs,amscd}
\usepackage{dsfont}
\usepackage{cite}
\usepackage{hyperref}
\usepackage{authblk}

\newtheorem{T}{Theorem}

\theoremstyle{remark}
\newtheorem{remark}{Remark}

\newcommand*{\cm}[1]{\mathscr{#1}}

\providecommand{\keywords}[1]
{
  \small	
  \textit{Key Words:~} #1
}

\providecommand{\MSC}[1]
{
  \small	
  \textit{Mathematics Subject Classification 2020:~} #1
}

\begin{document}

\title{Transition energy fields in the method of correlation equations}

\author{Linda A. Khachatryan\footnote{corresponding author, e-mail: \href{mailto: linda@instmath.sci.am}{linda@instmath.sci.am}}, Boris S. Nahapetian}

\affil{\small Institute of Mathematics, National Academy of Science of the Republic of Armenia}

\date{}

\maketitle

\begin{abstract}
In this paper, the well-known method of correlation equations for constructing Gibbs measures is generalized based on the concept of the transition energy field. Using the properties of transition energies, we obtain the system of correlation equations for lattice systems with finite spin space. It is shown that for a sufficiently small value of the one-point transition energies, the corresponding system of correlation functions, considered in infinite space, has a solution which is unique. Finally, the convergence of finite-volume correlation functions to the limiting correlation function is shown.
\end{abstract}


\MSC{82B05, 82B20, 60G60}

\section*{Introduction}

The method of correlation equations is one of the most demanded means of constructing and studying systems of mathematical statistical physics in infinite volumes (see, for example,~\cite{GMS, M, M2, R1, R2}). In the case of lattice systems, this method has generally been applied to vacuum systems which consists of two points --- spin and vacuum.

The problem of extending the method of correlation equations for more general systems naturally arises. This problem was considered in~\cite{N75}, in which a measurable set of finite measure was considered as a spin space, and the vacuum measure was taken to be equal to unity.

In all mentioned papers, the definition of correlation function is based on the notion of interaction potential. First, the case of pair-interaction potentials was considered (see, for example,~\cite{R1, M}). Further generalization was made by considering many body potentials (see~\cite{GMS, N75}).

In~\cite{DN19}, the new point of view on the mathematical foundations of statistical physics of infinite volume systems was presented. This viewpoint was based on the introduced notion of transition energy field. Applications of this notion in the theory of random fields was further developed in~\cite{KhN22, KhN23, KhN25}. In particular, it was shown that Gibbsian representation in terms of transition energies is universal for any finite-volume probability distributions.

There are various methods for constructing transition energy fields (see~\cite{DN19} and~\cite{KhN25}) including construction based on the axiomatic definition of Hamiltonian introduced in~\cite{DN19}, which does not assume that Hamiltonian is represented as a sum of interaction potential. This fact makes it possible to generalize known results of mathematical statistical physics.

In the present paper, we consider systems with finite spin space. Based on the results of~\cite{DN19, KhN22, KhN23, KhN25}, a system of correlation equations is written using the concept of the transition energy field. It is shown that for a sufficiently small value of the one-point transition energies, the corresponding system of correlation functions, considered in infinite space, has a solution which is unique. Finally, we prove the convergence of finite-volume correlation functions to limiting correlation function. For the particular case corresponding to finite-range interaction, these questions were considered in~\cite{KhN25-Reports}.

Note also that our approach can be applied to the systems with spin space considered in~\cite{N75}.

\section{Notations}

Let $\mathbb{Z}^d$ be a $d$-dimensional integer lattice, i.e., a set of $d$-dimensional vectors with integer components, $d \geq 1$. Note that all the arguments in this paper remain valid if we consider an arbitrary countable set instead of $\mathbb{Z}^d$.

For $S \subset \mathbb{Z}^d$, denote by $W(S) = \{ V \subset S , 0 < \left| V \right| < \infty \}$ the set of all nonempty finite subsets of $S$, where $\left| \Lambda \right|$ is the number of points in $\Lambda$. In the case $S = \mathbb{Z}^d$, we will use the simpler notation $W$. To denote the complement of the set $S$, we will write $S^c$. For one-point sets $\{t\}$, $t \in \Lambda$, braces will be omitted.

For $t = (t^{(1)},t^{(2)},...t^{(d)}), s = (s^{(1)},s^{(2)},...s^{(d)}) \in \mathbb{Z}^d$, we denote $\vert t-s \vert = \max \limits_{1 \le j \le d} \vert t^{(j)} - s^{(j)} \vert$. Also, we define $d(T,S) = \inf \limits_{t \in T, s \in S} \vert t-s \vert$, $T, S \subset \mathbb{Z}^d$.

Let each point $t \in \mathbb{Z}^d$ be associated with a set $X^t$, which is a copy of some finite set $X$, $1< \vert X \vert <\infty$. Denote by $X^S$ the set of all configurations on $S$, $S \subset \mathbb{Z}^d$, that is, the set $X^S = \{ x = ( x_t ,t \in S ), x_t \in X \}$, of all functions defined on $S$ and tacking values in $X$. For $S = {\O}$, we assume that $X^{\O} = \{ \boldsymbol{{\O}} \}$ where $\boldsymbol{{\O}}$ is an empty configuration. For any disjoint $S, T \subset \mathbb{Z}^d$ and any $x \in X^S$, $y \in X^T$, denote by $xy$ the concatenation of $x$ and $y$, that is, the configuration on $S \cup T$ equal to $x$ on $S$ and to $y$ on $T$. When $T \subset S$, we denote by $x_T$ the restriction of configuration $x \in X^S$ on $T$, i.e., $x_T=(x_t, t \in T)$.

Let $\theta_t$ be some fixed element of $X^t$ (vacuum) and $\theta = (\theta_t, t \in \mathbb{Z}^d)$. Denote $X_*^t = X^t \backslash \{\theta_t\}$, $t \in \mathbb{Z}^d$. For any $S \subset \mathbb{Z}^d$, denote by $X^S_*$ the set of configurations in $S$ which components do not contain the vacuum, and let $L^S_* = \bigcup \limits_{J \in W(S)} X^J_*$ be the set of configurations without vacuum which supports are subsets of~$S$. In the case $S = \mathbb{Z}^d$, we denote $L_* =L^{\mathbb{Z}^d}_*$. Note that any configuration from $X^S$ can be written as $x \theta_{S \backslash I}$ where $x \in X^I_*$, $I \subset S$.

\section{Transition energy fields}

In~\cite{DN19}, the notions of transition energy field and one-point transition energy field were introduced.

A set $\Delta = \{ \Delta_\Lambda^{\bar x}, \bar x \in X^{\Lambda^c}, \Lambda \in W \}$ of functions $\Delta_\Lambda^{\bar x}(x,u)$, $x,u \in X^\Lambda$, is called \emph{transition energy field} if its elements satisfy the following consistency conditions: for all $\Lambda \in W$ and $\bar x \in X^{\Lambda^c}$, it holds
$$
\Delta_\Lambda^{\bar x}(x,u) = \Delta_\Lambda^{\bar x}(x,y) + \Delta_\Lambda^{\bar x}(y,u), \qquad x,u,y \in X^\Lambda;
$$
and for all disjoint $\Lambda, V \in W$ and $\bar x \in X^{(\Lambda \cup V)^c}$,
$$
\Delta_{\Lambda \cup V}^{\bar x}(xy,uv) = \Delta_\Lambda^{\bar x y}(x,u) + \Delta_V^{\bar x u}(y,v), \qquad x,u \in X^\Lambda, y,v \in X^V.
$$

Note that in particular, it holds
$$
\Delta_\Lambda^{\bar x}(x,u) = - \Delta_\Lambda^{\bar x}(u,x), \qquad \Delta_\Lambda^{\bar x}(x,x) = 0, \qquad x,u \in X^\Lambda,
$$
and
$$
\Delta_{\Lambda \cup V}^{\bar x}(xy,uy) = \Delta_\Lambda^{\bar x y}(x,u) \qquad x,u \in X^\Lambda, y \in X^V.
$$

A set $\Delta_1 = \{ \Delta_t^{\bar x}, \bar x \in X^{t^c}, t \in \mathbb{Z}^d \}$ of functions $\Delta_t^{\bar x}(x,u)$, $x,u \in X^t$, is called \emph{one-point transition energy field} if its elements satisfy the following consistency conditions: for all $t \in \mathbb{Z}^d$ and $\bar x \in X^{t^c}$, it holds
$$
\Delta_t^{\bar x}(x,u) = \Delta_t^{\bar x}(x,y) + \Delta_t^{\bar x}(y,u), \qquad x,u,y \in X^t;
$$
and for all $t,s \in \mathbb{Z}^d$ and $\bar x \in X^{\{t,s\}^c}$,
$$
\Delta_t^{\bar x y}(x,u) + \Delta_s^{\bar x u}(y,v) = \Delta_s^{\bar x x}(y,v) + \Delta_t^{\bar x v}(x,u), \qquad x,u \in X^t, y,v \in X^s.
$$

The following theorem states the relationship between the elements of a transition energy field and a one-point transition energy field (see~\cite{DN19} as well as~\cite{KhN23}).

\begin{T}
\label{th-Delta-delta}
A set $\Delta = \{ \Delta_\Lambda^{\bar x}, \bar x \in X^{\Lambda^c}, \Lambda \in W \}$ of functions on $X^\Lambda \times X^\Lambda$, $\Lambda \in W$, is a transition energy field if and only if its elements can be represented in the form
\begin{equation}
\label{Delta-Delta1}
\Delta_\Lambda^{\bar x}(x,u) = \Delta_{t_1}^{\bar x x_{\{\Lambda \backslash t_1\}}} (x_{t_1},u_{t_1}) + \Delta_{t_2}^{\bar x u_{t_1} x_{\{\Lambda \backslash \{t_1,t_2\}\}}} (x_{t_2},u_{t_2}) + \ldots + \Delta_{t_n}^{\bar x u_{\{\Lambda \backslash t_n\}}} (x_{t_n},u_{t_n}),
\end{equation}
where $x,u \in X^\Lambda$, $\Lambda = \{t_1, t_2, ..., t_n\}$ is some enumeration of points in $\Lambda$, $\vert \Lambda \vert = n$, and $\Delta_1 = \{\Delta_t^{\bar x}, \bar x \in X^{t^c}, t \in \mathbb{Z}^d\}$ is a one-point transition energy field.
\end{T}

Thus, the one-point transition energy field $\Delta_1$ uniquely determines the transition energy field~$\Delta$. Therefore, when obtaining results, conditions can only be imposed on $\Delta_1$.

\section{Correlation functions}

Let $\Delta_1 = \{\Delta_t^{\bar x}, \bar x \in X^{t^c}, t \in \mathbb{Z}^d\}$ be a one-point transition energy field and $\Delta = \{ \Delta_\Lambda^{\bar x}, \bar x \in X^{\Lambda^c}, \Lambda \in W \}$ be the corresponding transition energy field. Let us fix some $\Lambda \in W$. To simplify notations, we denote $\Delta_\Lambda = \Delta_\Lambda^{\theta_{\Lambda^c}}$, and for any $t \in \Lambda$ and $z \in X^I$, $I \subset \Lambda \backslash t$, we will write $\Delta_t^z$ instead of $\Delta_t^{z \theta_{(t \cup I)^c}}$.

For the given $\Delta_1$, we define finite-volume correlation function relative to $\Lambda$ as a function $\rho_\Lambda$ on $L_*$ by the formula
$$
\rho_\Lambda(x) = \left\{
  \begin{array}{ll}
    Z_\Lambda^{-1} \sum \limits_{y \in X^{\Lambda \backslash I}} \exp\{\Delta_\Lambda(x y, \theta_\Lambda)\}, & x \in X_*^I, I \subset \Lambda, \\
    0, & x \in X_*^I, I \not\subset \Lambda, \\
  \end{array}     \right.
$$
where
$$
Z_\Lambda = \sum \limits_{x \in X^\Lambda} \exp\{\Delta_\Lambda(x, \theta_\Lambda)\},
$$
and $\rho_\Lambda (\boldsymbol{\O}) =1$.

Thus, each $\Delta_1$ defines a set $\rho(\Delta_1) = \{\rho_\Lambda, \Lambda \in W\}$ of finite-volume correlation functions. In the theorem below, we show that under a suitable condition on the elements of $\Delta_1$, each finite-volume correlation function satisfies a certain equation. This condition imposes a restriction on the boundary conditions (environment) in the following sense: for any $\Lambda \in W$, $t,s \in \Lambda$ and $I \subset \Lambda \backslash \{t,s\}$, it holds
\begin{equation}
\label{Delta-Phi-2}
\Delta_t^{zy}(x, \theta_t) - \Delta_t^{zv}(x, \theta_t) = \Delta_t^y(x, \theta_t) - \Delta_t^v(x, \theta_t),
\end{equation}
where $x \in X^t$, $y,v \in X^s$, $z \in X^I$.

Finally, we note that the proof of the following theorem is based solely on the properties of transition energies and does not assume that the energy is the sum of interactions.

\begin{T}
\label{prop-rho-Lambda}
Let $\Delta_1 = \{\Delta_t^{\bar x}, \bar x \in X^{t^c}, t \in \mathbb{Z}^d\}$ be a one-point transition energy field satisfying~\eqref{Delta-Phi-2}. Then for any $\Lambda \in W$, correlation function $\rho_\Lambda$ satisfies the following equation: for any $t \in I \subset \Lambda \in W$ and $x \in X_*^t$, $u \in X_*^{I \backslash t}$, it holds
$$
\rho_\Lambda (xu) = \dfrac{e^{\Delta_t^u(x,\theta_t)}}{\sum \limits_{\alpha \in X^t} e^{\Delta_t^u(\alpha,\theta_t)}} \left( \rho_\Lambda (u) + \sum \limits_{\alpha \in X^t} e^{\Delta_t^u(\alpha,\theta_t)} \Big( G_\Lambda(xu) - G_\Lambda (\alpha u) \Big) \right),
$$
where
$$
G_\Lambda(xu) = \sum \limits_{J \subset \Lambda \backslash I} \sum \limits_{y \in X_*^J} K_{t \cup J}(xy) \left( \rho_\Lambda(uy) - \sum \limits_{\alpha \in X_*^t} \rho_\Lambda(\alpha uy) \right)
$$
and
$$
K_{t \cup J}(xy) = \prod \limits_{s \in J} \left( e^{\Delta_s^x(y_s,\theta_s) - \Delta_s(y_s,\theta_s)} - 1\right).
$$
\end{T}
\begin{proof}
Let $t \in I \subset \Lambda \in W$ and $x \in X_*^t$, $u \in X_*^{I \backslash t}$. To shorten notations, denote $\theta = \theta_\Lambda$. According to the properties of transition energies, for any $y \in X^{\Lambda \backslash I}$, we have
$$
\begin{array}{l}
\Delta_\Lambda(xuy, \theta) = \Delta_\Lambda(xuy, \theta_t uy) + \Delta_\Lambda(\theta_t uy, \theta) + \Delta_\Lambda(xu \theta_{\Lambda \backslash I}, \theta_t u \theta_{\Lambda \backslash I}) - \Delta_\Lambda(xu \theta_{\Lambda \backslash I}, \theta_t u \theta_{\Lambda \backslash I})=\\
\\
= \Delta_t^{uy}(x, \theta_t) + \Delta_\Lambda(\theta_t uy, \theta) + \Delta_t^u(x,\theta_t) - \Delta_t^u(x,\theta_t).
\end{array}
$$
Therefore, we can write
$$
\begin{array}{l}
\rho_\Lambda (xu) = \dfrac{1}{Z_\Lambda} \sum \limits_{y \in X^{\Lambda \backslash I}} \exp\{\Delta_\Lambda(xuy, \theta)\} = \\
\\
= \dfrac{e^{\Delta_t^u(x,\theta_t)}}{Z_\Lambda} \sum \limits_{y \in X^{\Lambda \backslash I}} e^{\Delta_\Lambda(\theta_t uy, \theta)} \left( 1 + \exp\{\Delta_t^{uy}(x,\theta_t) - {\Delta_t^{u}(x,\theta_t)}\} - 1  \right) = \\
\\
= e^{\Delta_t^u(x,\theta_t)} \left(\dfrac{1}{Z_\Lambda} \sum \limits_{y \in X^{\Lambda \backslash I}} e^{\Delta_\Lambda(\theta_t uy, \theta)} + \dfrac{1}{Z_\Lambda} \sum \limits_{y \in X^{\Lambda \backslash I}} \left( e^{\Delta_t^{uy}(x,\theta_t) - \Delta_t^{u}(x,\theta_t)} - 1  \right) e^{\Delta_\Lambda(\theta_t uy, \theta)} \right).
\end{array}
$$
Since
$$
e^{\Delta_\Lambda(\theta_t uy, \theta)} = \sum \limits_{\alpha \in X^t} e^{\Delta_\Lambda(\alpha uy, \theta)} - \sum \limits_{\alpha \in X_*^t} e^{\Delta_\Lambda(\alpha uy, \theta)},
$$
for the first summand in the obtained relation, we have
$$
\dfrac{1}{Z_\Lambda} \sum \limits_{y \in X^{\Lambda \backslash I}} e^{\Delta_\Lambda(\theta_t uy, \theta)} = \rho_\Lambda (u) - \sum \limits_{\alpha \in X_*^t} \rho_\Lambda (\alpha u).
$$

Let us consider the second summand
$$
\begin{array}{l}
G_\Lambda(xu) = \dfrac{1}{Z_\Lambda} \sum \limits_{y \in X^{\Lambda \backslash I}} \left( e^{\Delta_t^{uy}(x,\theta_t) - \Delta_t^{u}(x,\theta_t)} - 1 \right) e^{\Delta_\Lambda(\theta_t uy, \theta)} = \\
\\
= \dfrac{1}{Z_\Lambda} \sum \limits_{J \subset \Lambda \backslash I} \sum \limits_{y \in X_*^J} \left( e^{\Delta_t^{uy}(x,\theta_t) - \Delta_t^{u}(x,\theta_t)} - 1 \right) e^{\Delta_\Lambda(uy \theta_{t \cup (\Lambda \backslash I \backslash J)}, \theta)}.
\end{array}
$$
Using properties of $\boldsymbol{\Delta}$, we can write
$$
\begin{array}{l}
\Delta_t^{uy}(x,\theta_t) - \Delta_t^{u}(x,\theta_t) = \Delta_\Lambda(xuy \theta_{\Lambda \backslash I \backslash J}, uy \theta_{t \cup (\Lambda \backslash I \backslash J)}) - \Delta_\Lambda(xu \theta_{\Lambda \backslash I}, u \theta_{t \cup (\Lambda \backslash I)}) = \\
\\
= \Delta_\Lambda(xuy \theta_{\Lambda \backslash I \backslash J}, xu \theta_{\Lambda \backslash I}) + \Delta_\Lambda(xu \theta_{\Lambda \backslash I}, u \theta_{t \cup (\Lambda \backslash I)}) + \Delta_\Lambda(u \theta_{t \cup (\Lambda \backslash I)}, uy \theta_{t \cup (\Lambda \backslash I \backslash J)}) - \\
\\
- \Delta_\Lambda(xu \theta_{\Lambda \backslash I}, u \theta_{t \cup (\Lambda \backslash I)}) = \\
\\
= \Delta_J^{xu}(y,\theta_J) - \Delta_J^u(y,\theta_J) = \sum \limits_{j = 1}^{\vert J \vert} \left( \Delta_{s_j}^{xu y_{s_{j+1},...,s_{\vert J \vert}}}(y_{s_j}, \theta_{s_j}) - \Delta_{s_j}^{u y_{s_{j+1},...,s_{\vert J \vert}}}(y_{s_j}, \theta_{s_j})\right),
\end{array}
$$
where we used formula~\eqref{Delta-Delta1} with $J = \{s_1, s_2,...,s_{\vert J \vert}\}$ being some enumeration of points in $J$. Tacking into account condition~\eqref{Delta-Phi-2}, we obtain
$$
\Delta_t^{uy}(x,\theta_t) - \Delta_t^{u}(x,\theta_t) = \sum \limits_{s \in J} \left( \Delta_s^x(y_s, \theta_s) - \Delta_s(y_s, \theta_s) \right).
$$
Hence,
$$
\begin{array}{l}
G_\Lambda(xu) = \dfrac{1}{Z_\Lambda} \sum \limits_{J \subset \Lambda \backslash I} \sum \limits_{y \in X_*^J} \left( \prod \limits_{s \in J} e^{\Delta_s^x(y_s, \theta_s) - \Delta_s(y_s, \theta_s)} - 1 \right) e^{\Delta_\Lambda(uy \theta_{t \cup (\Lambda \backslash I \backslash J)}, \theta)} = \\
\\
= \dfrac{1}{Z_\Lambda} \sum \limits_{J \subset \Lambda \backslash I} \sum \limits_{y \in X_*^J} \sum \limits_{T \subset J} \prod \limits_{s \in T} \left(e^{\Delta_s^x(y_s, \theta_s) - \Delta_s(y_s, \theta_s)} - 1 \right) e^{\Delta_\Lambda(uy \theta_{t \cup (\Lambda \backslash I \backslash J)}, \theta)} = \\
\\
= \dfrac{1}{Z_\Lambda} \sum \limits_{T \subset \Lambda \backslash I} \sum \limits_{J \subset \Lambda \backslash I: T \subset J} \sum \limits_{y \in X_*^T} \sum \limits_{z \in X_*^J} \prod \limits_{s \in T} \left(e^{\Delta_s^x(y_s, \theta_s) - \Delta_s(y_s, \theta_s)} - 1 \right) e^{\Delta_\Lambda(uyz \theta_{t \cup (\Lambda \backslash I \backslash J)}, \theta)} =\\
\\
= \sum \limits_{T \subset \Lambda \backslash I} \sum \limits_{y \in X_*^T} \prod \limits_{s \in T} \left(e^{\Delta_s^x(y_s, \theta_s) - \Delta_s(y_s, \theta_s)} - 1 \right) \dfrac{1}{Z_\Lambda} \sum \limits_{J \subset \Lambda \backslash I \backslash T} \sum \limits_{z \in X_*^J} e^{\Delta_\Lambda(uyz \theta_{t \cup (\Lambda \backslash I \backslash J \backslash T)}, \theta)}.
\end{array}
$$
Since
$$
\dfrac{1}{Z_\Lambda} \sum \limits_{J \subset \Lambda \backslash I \backslash T} \sum \limits_{z \in X_*^J} e^{\Delta_\Lambda(uyz \theta_{t \cup (\Lambda \backslash I \backslash J \backslash T)}, \theta)} = \dfrac{1}{Z_\Lambda} \sum \limits_{z \in X^{\Lambda \backslash I \backslash T}} e^{\Delta_\Lambda(\theta_t uyz, \theta)} = \rho_\Lambda (uy) - \sum \limits_{\alpha \in X_*^t} \rho_\Lambda (\alpha uy),
$$
finally we obtain
$$
G_\Lambda(xu) = \sum \limits_{J \subset \Lambda \backslash I} \sum \limits_{y \in X_*^J} K_{t \cup J}(xy) \left( \rho_\Lambda(uy) - \sum \limits_{\alpha \in X_*^t} \rho_\Lambda(\alpha uy) \right)
$$
where
$$
K_{t \cup J}(xy) = \prod \limits_{s \in J} \left( e^{\Delta_s^x(y_s,\theta_s) - \Delta_s(y_s,\theta_s)} - 1\right).
$$

Thus,
$$
\rho_\Lambda (xu) = e^{\Delta_t^u(x,\theta_t)} \left( \rho_\Lambda (u) - \sum \limits_{\alpha \in X_*^t} \rho_\Lambda (\alpha u) + G_\Lambda(xu) \right).
$$
Tacking the sum of the both sides of the obtained relation over all $x \in X_*^t$, we get
$$
\sum \limits_{x \in X_*^t} \rho_\Lambda (xu) = \rho_\Lambda (u) \sum \limits_{x \in X_*^t} e^{\Delta_t^u(x,\theta_t)} - \left( \sum \limits_{x \in X_*^t} e^{\Delta_t^u(x,\theta_t)} \right) \cdot \sum \limits_{\alpha \in X_*^t} \rho_\Lambda (\alpha u) + \sum \limits_{x \in X_*^t} e^{\Delta_t^u(x,\theta_t)} G_\Lambda (xu),
$$
and therefore,
$$
\sum \limits_{\alpha \in X_*^t} \rho_\Lambda (\alpha u) = \dfrac{\sum \limits_{\alpha \in X_*^t} e^{\Delta_t^u(\alpha,\theta_t)}}{1 + \sum \limits_{\alpha \in X_*^t} e^{\Delta_t^u(\alpha,\theta_t)}} \rho_\Lambda (u) + \dfrac{1}{1 + \sum \limits_{\alpha \in X_*^t} e^{\Delta_t^u(\alpha,\theta_t)}} \sum \limits_{\alpha \in X_*^t}  e^{\Delta_t^u(\alpha,\theta_t)} G_\Lambda (\alpha u).
$$
Since $K_{t \cup J}(\theta_t y) = 0$, we have $G_\Lambda(\theta_t u) = 0$, and, finally, we obtain
$$
\begin{array}{l}
\rho_\Lambda (xu) = e^{\Delta_t^u(x,\theta_t)} \left( \dfrac{1}{\sum \limits_{\alpha \in X^t} e^{\Delta_t^u(\alpha,\theta_t)}} \rho_\Lambda (u) - \dfrac{1}{\sum \limits_{\alpha \in X^t} e^{\Delta_t^u(\alpha,\theta_t)}} \sum \limits_{\alpha \in X_*^t} e^{\Delta_t^u(\alpha,\theta_t)} G_\Lambda (\alpha u) + G_\Lambda(xu) \right) = \\
\\
= \dfrac{e^{\Delta_t^u(x,\theta_t)}}{\sum \limits_{\alpha \in X^t} e^{\Delta_t^u(\alpha,\theta_t)}} \left( \rho_\Lambda (u) + \sum \limits_{\alpha \in X^t} e^{\Delta_t^u(\alpha,\theta_t)} \Big( G_\Lambda(xu) - G_\Lambda (\alpha u) \Big) \right).
\end{array}
$$
\end{proof}

\section{Equations for correlation functions}

Let $\Delta_1 = \{\Delta_t^{\bar x}, \bar x \in X^{t^c}, t \in \mathbb{Z}^d\}$ be a one-point transition energy field. We introduce the norm for $\Delta_1$ as follows:
$$
\|\Delta_1\| = \sup \limits_{t \in \mathbb{Z}^d} \sup \limits_{x,u \in X^t} \sup \limits_{\bar x \in X^{t^c}} \vert \Delta_t^{\bar x} (x, u) \vert.
$$
Put also
$$
D = \sup \limits_{t \in \mathbb{Z}^d} \sup \limits_{x \in X^t} \sum \limits_{s \in t^c} \sup \limits_{y \in X^s} \vert \Delta_s^x (y, \theta_s) - \Delta_s (y, \theta_s) \vert.
$$

Consider the Banach space $\cm{B}_*$ of bounded functions $\varphi$ on $L_*$ with the norm
$$
\| \varphi \| = \sup_{\Lambda \in W} \| \varphi \|_\Lambda, \qquad
\| \varphi \|_\Lambda = \sum \limits_{x \in X_*^\Lambda} \vert \varphi (x) \vert.
$$

We assume that $\mathbb{Z}^d$ is endowed with some order $\preceq$, for example, the lexicographical order. For each $I \in W$, denote $I' = I \backslash t$ where $t$ is the smallest element in $I$ with respect to $\preceq$. For the sake of simplicity, for any $x \in X_*^I$, we will use the notation $x' = x_{I'}$.

Consider the operator $\cm{K} = \cm{K}^{\Delta_1}$ on $\cm{B}_*$ defined as follows:
$$
(\cm{K} \varphi)(x) = \gamma(x) \big((S \varphi)(x) + (T \varphi)(x)\big), \qquad x \in X_*^I, I \in W,
$$
where
$$
\gamma(x) = \dfrac{e^{\Delta_t^{x'}(x_t,\theta_t)}}{\sum \limits_{\alpha \in X^t} e^{\Delta_t^{x'}(\alpha,\theta_t)}},
\qquad
(S \varphi)(x) = \left\{ \begin{array}{l}
                              \varphi(x'), \quad \vert I \vert > 1, \\
                              \\
                              0, \qquad \quad \vert I \vert = 1,
                            \end{array} \right.
$$
and
$$
(T \varphi)(x) = (G \varphi)(x) + \sum \limits_{\alpha \in X_*^t} e^{\Delta_t^{x'}(\alpha,\theta_t)} \Big( (G \varphi)(x) - (G \varphi) (\alpha x') \Big)
$$
with
$$
(G \varphi)(x) = \sum \limits_{J \in W(I^c)} \sum \limits_{y \in X_*^J} K_{t \cup J}(x_t y) \left( \varphi(x'y) - \sum \limits_{\alpha \in X_*^t} \varphi(\alpha x'y) \right)
$$
and
$$
K_{t \cup J}(x_t y) = \prod \limits_{s \in J} \left( e^{\Delta_s^{x_t} (y_s,\theta_s) - \Delta_s(y_s,\theta_s)} - 1\right).
$$

Further, we put
$$
\delta (x) = \left\{ \begin{array}{ll}
                       \gamma(x), & \vert I \vert = 1, \\
                       0, & \vert I \vert > 1,
                     \end{array}
 \right. \qquad x \in X^I_*, I \in W,
$$
and
$$
C_1 = \frac{e^{\| \Delta_1 \|} N_X}{1 + e^{\| \Delta_1 \|} N_X}, \qquad C_2 = 2 \big( 1 + 2 e^{\| \Delta_1 \|} N_X \big) \big(\exp \{e^D -1\} -1 \big), \qquad N_X = \vert X \vert - 1.
$$

\begin{T}
\label{prop-K-norm}
Let $\Delta_1$ be a one-point transition energy field satisfying~\eqref{Delta-Phi-2} and such that
\begin{equation}
\label{Delta<1}
C_1 (1 + C_2) < 1.
\end{equation}
Then the equation
\begin{equation}
\label{cor-eq}
\rho = \delta + \cm{K} \rho
\end{equation}
with $\cm{K} = \cm{K}^{\Delta_1}$ has a unique solution on $\cm{B}_*$ given by
$$
\rho = \big( 1-  \cm{K} \big)^{-1} \delta = \left( 1 + \sum \limits_{n=1}^\infty \cm{K}^n \right) \delta.
$$
\end{T}
\begin{proof}
The mentioned solution can be obtained if we subsequently apply operator $\cm{K}$ to both sides of equation~\eqref{cor-eq}. It remains to show that $\|\cm{K}\|<1$.

Let $I \in W$ and $x \in X^I_*$. Since
$$
- \| \Delta_1 \| \le \Delta_t^{x'}(x_t,\theta_t) \le \| \Delta_1 \|,
$$
we have
$$
\vert \gamma (x) \vert = \left| \dfrac{e^{\Delta_t^{x'}(x_t,\theta_t)}}{1 + \sum \limits_{\alpha \in X_*^t} e^{\Delta_t^{x'}(\alpha,\theta_t)}} \right| \le \frac{e^{\| \Delta_1 \|}}{1 + e^{-{\| \Delta_1 \|}} N_X}.
$$
Hence, $\| \gamma S \varphi \|_I = 0$ if $\vert I \vert = 1$, and for $\vert I \vert > 1$, we have
$$
\| \gamma S \varphi \|_I = \sum \limits_{x \in X_*^I} \vert \gamma (x) \vert \cdot \vert (S\varphi)(x) \vert \le \frac{e^{\| \Delta_1 \|}}{1 + e^{-{\| \Delta_1 \|}} N_X} \sum \limits_{x \in X^t_*}  \sum \limits_{x' \in X_*^{I'}} \vert \varphi(x') \vert = C_1 \|\varphi\|_{I'}.
$$
Thus, $\| \gamma S \varphi \|_I \le C_1 \|\varphi\|$ for any $I \in W$.

Further, we can write
$$
\begin{array}{l}
\|\gamma T \varphi\|_I = \sum \limits_{x \in X_*^I} \vert \gamma(x) \vert \left| (G \varphi)(x) \left( 1 + \sum \limits_{\alpha \in X_*^t} e^{\Delta_t^{x'}(\alpha,\theta_t)} \right) -  \sum \limits_{\alpha \in X_*^t} e^{\Delta_t^{x'}(\alpha,\theta_t)} (G \varphi) (\alpha x') \right| \le \\
\\
\le C_1 \big( 1 + 2 e^{\| \Delta_1 \|} N_X \big) \sum \limits_{x' \in X_*^{I'}} \left| \sup \limits_{\alpha \in X_*^t} (G \varphi)(\alpha x') \right| \le \\
\\
\le C \sum \limits_{x' \in X_*^{I'}} \sum \limits_{J \in W(I^c)} \sum \limits_{y \in X_*^J} \left| \sup \limits_{\alpha \in X_*^t} K_{t \cup J}(\alpha y) \right| \cdot \left| \varphi (x'y) - \sum \limits_{\alpha \in X_*^t} \varphi (\alpha x'y) \right| \le \\
\\
= C \sum \limits_{J \in W(I^c)} \sup \limits_{\alpha \in X_*^t, y \in X_*^J} \vert K_{t \cup J} (\alpha y) \vert \left( \|\varphi\|_{I' \cup J}  + \|\varphi\|_{I \cup J} \right) \le \\
\\
\le 2 \| \varphi \| C  \sum \limits_{J \in W(I^c)} \sup \limits_{\alpha \in X_*^t, y \in X_*^J} \vert K_{t \cup J} (\alpha y) \vert,
\end{array}
$$
where $C = C_1 \big( 1 + 2 e^{\| \Delta_1 \|} N_X \big)$. Since for any $\alpha \in X_*^t$ and $y \in X_*^{I^c}$, using properties of exponential function, we can write
$$
\begin{array}{l}
\sum \limits_{J \in W(I^c)} \vert K_{t \cup J} (\alpha y_J) \vert = \sum \limits_{J \in W(I^c)} \prod \limits_{s \in J} \vert  e^{\Delta_s^\alpha (y_s,\theta_s) - \Delta_s(y_s,\theta_s)} - 1 \vert = \\
\\
= \sum \limits_{n=1}^{\infty} \sum \limits_{J \in W(I^c) : \vert J \vert = n} \prod \limits_{s \in J} \vert  e^{\Delta_s^\alpha (y_s,\theta_s) - \Delta_s(y_s,\theta_s)} - 1 \vert = \\
\\
= \sum \limits_{n=1}^{\infty} \sum \limits_{s_1 \in I^c} \sum \limits_{s_2 \in I^c \backslash s_1} ... \sum \limits_{s_n \in I^c \backslash \{s_1, s_2, ..., s_{n-1}\}} \prod \limits_{j=1}^n \vert  e^{\Delta_{s_j}^\alpha (y_{s_j},\theta_{s_j}) - \Delta_{s_j} (y_{s_j},\theta_{s_j})} - 1 \vert < \\
\\
< \sum \limits_{n=1}^{\infty} \dfrac{1}{n!} \left( \sum \limits_{s \in I^c} \vert  e^{\Delta_s^\alpha (y_s,\theta_s) - \Delta_s (y_s,\theta_s)} - 1 \vert \right)^n \le \sum \limits_{n=1}^{\infty} \dfrac{1}{n!} \left(  e^{\sum \limits_{s \in I^c} \vert \Delta_s^\alpha (y_s,\theta_s) - \Delta_s (y_s,\theta_s)\vert} - 1 \right)^n < \\
\\
< \sum \limits_{n=1}^\infty \dfrac{1}{n!} \left(  e^D - 1 \right)^n = \exp \{e^D -1\} -1,
\end{array}
$$
finally we obtain
$$
\|\gamma T \varphi\|_I \le 2 C_1 \big( 1 + 2 e^{\| \Delta_1 \|} N_X \big) \big(\exp \{e^D -1\} -1 \big) \| \varphi \|
$$
for any $I \in W$.

Therefore,
$$
\| \cm{K}\| \le C_1 + 2 C_1 \big( 1 + 2 e^{\| \Delta_1 \|} N_X \big) \big(\exp \{e^D -1\} -1 \big) = C_1 (1 + C_2) < 1.
$$
\end{proof}

The following theorem is the main result of the paper. In its proof, we follow Ruelle~\cite{R1, R2}.

\begin{T}
\label{th-main}
Let $\Delta_1$ be a one-point transition energy field satisfying~\eqref{Delta-Phi-2} and~\eqref{Delta<1} and let $\rho(\Delta_1) = \{\rho_\Lambda, \Lambda \in W\}$ be the set of corresponding finite-volume correlation functions. Then there exists a positive decreasing function $\varepsilon$ such that
$$
\lim \limits_{d \to \infty} \varepsilon(d) = 0
$$
and for any $x \in X^I_*$, $I \subset \Lambda \in W$,
$$
\left| \rho_\Lambda(x) - \rho(x) \right| < \varepsilon \big( d(I, \Lambda^c) \big),
$$
where $\rho$ is the solution to equation~\eqref{cor-eq}.
\end{T}
\begin{proof}
For any $\Lambda \in W$, consider the operator $\psi_\Lambda$ on $\cm{B}_*$ defined as follows:
$$
(\psi_\Lambda \varphi)(x) = \left\{ \begin{array}{lr}
                                      \varphi(x), & x \in L_*^\Lambda, \\
                                      0, & \text{otherwise}.
                                    \end{array} \right.
$$
It is clear that $\|\psi_\Lambda\| = \sup \limits_{\|\varphi\| = 1} \|\psi_\Lambda \varphi\| = 1$, and for any $V \subset \Lambda$, one has $\psi_V \psi_\Lambda = \psi_V$.

According to Theorem~\ref{prop-rho-Lambda}, for each $\Lambda \in W$, the corresponding element $\rho_\Lambda$ of $\rho(\Delta_1)$ satisfies the equation
\begin{equation}
\label{cor-eq-Lambda}
\rho_\Lambda = \psi_\Lambda \delta + \psi_\Lambda \cm{K} \rho_\Lambda,
\end{equation}
Since $\|\psi_\Lambda \cm{K}\| \le \|\cm{K}\|$, due to Theorem~\ref{prop-K-norm} we have $\|\psi_\Lambda \cm{K}\| < 1$, and hence, $\rho_\Lambda$ is the unique solution to~\eqref{cor-eq-Lambda}, which can be written as
$$
\rho_\Lambda = \big( 1-\psi_\Lambda \cm{K} \big)^{-1} \psi_\Lambda \delta = \left( 1 + \sum \limits_{n=1}^\infty \big(\psi_\Lambda \cm{K} \big)^n \right) \psi_\Lambda \delta.
$$

For any $V \subset \Lambda \in W$ and $n \ge 1$, we have
$$
\begin{array}{l}
\psi_V \rho_\Lambda - \psi_V \rho = \psi_V \left( \big( 1-\psi_\Lambda \cm{K} \big)^{-1} \psi_\Lambda - \big( 1- \cm{K} \big)^{-1} \right) \delta = \\
\\
= \psi_V \left( \big( 1-\psi_\Lambda \cm{K} \big)^{-1} - \big( 1 + \sum \limits_{k=1}^n (\psi_\Lambda \cm{K} )^k \big) \right) \psi_\Lambda \delta + \psi_V \psi_\Lambda \delta + \psi_V \big( \sum \limits_{k=1}^n (\psi_\Lambda \cm{K} )^k \big) \psi_\Lambda \delta  - \\
\\
- \psi_V \left( \big( 1- \cm{K} \big)^{-1} - \big( 1 + \sum \limits_{k=1}^n \cm{K}^k \big) \right) \delta - \psi_V \delta - \psi_V \big( \sum \limits_{k=1}^n \cm{K}^k \big) \delta  = \\
\\
= \psi_V \big( \sum \limits_{k=n+1}^\infty (\psi_\Lambda \cm{K} )^k \big) \psi_\Lambda \delta - \psi_V \big( \sum \limits_{k=n+1}^\infty \cm{K}^k \big) \delta + \psi_V \sum \limits_{k=1}^n \big( (\psi_\Lambda \cm{K} )^k \psi_\Lambda - \cm{K}^k \big) \delta.
\end{array}
$$
Further,
$$
\left\| \psi_V \big( \sum \limits_{k=n+1}^\infty (\psi_\Lambda \cm{K} )^k \big) \psi_\Lambda \right\| \le \sum \limits_{k=n+1}^\infty \| \cm{K} \|^k = \frac{\| \cm{K} \|^{n+1}}{1-\| \cm{K} \|},
$$
$$
\left\| \psi_V \big( \sum \limits_{k=n+1}^\infty \cm{K}^k \big) \right\| \le \sum \limits_{k=n+1}^\infty \| \cm{K} \|^k = \frac{\| \cm{K} \|^{n+1}}{1-\| \cm{K} \|},
$$
and it remains to estimate $\left\| \psi_V \sum \limits_{k=1}^n \big( (\psi_\Lambda \cm{K} )^k \psi_\Lambda - \cm{K}^k \big) \right\|$.

For any $V \subset \Lambda \in W$ and $S \subset \mathbb{Z}^d$ such that $\Lambda \subset S$, and any $x \in X_*^I$, $I \subset V$, we can write
$$
\begin{array}{l}
\left| (\psi_V \cm{K} \psi_S \varphi)(x) - (\psi_V \cm{K} \psi_\Lambda \varphi)(x) \right| \le \\
\\
\le \left| \gamma(x)\right| \sum \limits_{\alpha \in X_*^t} e^{\Delta_t^{x'}(\alpha,\theta_t)}  \sum \limits_{ \substack {J \in W(S \backslash I):\\ J \cap (S \backslash \Lambda) \neq {\O} }} \sum \limits_{y \in X_*^J} \left| K_{t \cup J}(x_t y) - K_{t \cup J}(\alpha y) \right| \cdot \left| \varphi(x'y) - \sum \limits_{\beta \in X_*^t} \varphi(\beta x'y) \right| \le \\
\\
\le \dfrac{e^{\| \Delta_1 \|}}{1 + e^{-\| \Delta_1 \|} N_X} e^{\| \Delta_1 \|} N_X \sum \limits_{J \in W(\Lambda^c)} 2 \sup \limits_{\alpha \in X_*^t, y \in X_*^J}\left| K_{t \cup J}(\alpha y) \right| \left( \sum \limits_{y \in X_*^J} \vert \varphi(x'y)\vert + \sum \limits_{y \in X_*^J, \beta \in X_*^t} \vert \varphi(\beta x'y)\vert \right) \le \\
\\
\le 4 C_1 e^{\| \Delta_1 \|} \|\varphi\| \sum \limits_{J \in W(\Lambda^c)} \sup \limits_{\alpha \in X_*^t, y \in X_*^J}\left| K_{t \cup J}(\alpha y) \right|.
\end{array}
$$
Proceeding as in the proof of Theorem~\ref{prop-K-norm}, for any $\alpha \in X_*^t$ and $y \in X_*^{\Lambda^c}$, we obtain
$$
\sum \limits_{J \in W(\Lambda^c)} \left| K_{t \cup J}(\alpha y_J) \right| \le \sum \limits_{n=1}^{\infty} \dfrac{1}{n!} \left(  e^{\sum \limits_{s \in \Lambda^c} \vert \Delta_s^\alpha (y_s,\theta_s) - \Delta_s (y_s,\theta_s)\vert} - 1 \right)^n < \exp \{e^{\sigma(\Lambda)} -1\} -1,
$$
where
$$
\sigma(\Lambda) = \sum \limits_{s \in \Lambda^c} \sup \limits_{\alpha \in X^t, y \in X^s} \vert \Delta_s^\alpha (y,\theta_s) - \Delta_s (y,\theta_s)\vert \to 0 \quad \text{as } \Lambda \uparrow V^c.
$$
Therefore, there exists a positive decreasing function $f$ such that $\lim \limits_{d \to \infty} f(d) = 0$ and for any $V \subset \Lambda \subset S$,
\begin{equation}
\label{norm-K-dif}
\|\psi_V \cm{K} \psi_S - \psi_V \cm{K} \psi_\Lambda\| \le f(d(V,\Lambda^c)).
\end{equation}

For $\Lambda \in W$, denote by $\Lambda(r) = \{s \in \Lambda: d(s, \Lambda^c) > r\}$ the set of points of $\Lambda$ with distance greater than $r$ to $\Lambda^c$, $r \ge 1$. It is clear that for any $r_1 > r_2$, $\Lambda(r_1) \subset \Lambda (r_2) \subset \Lambda$. From~\eqref{norm-K-dif} it follows that for any $r \ge 1$ and $k \ge 1$, we have
$$
\| \psi_{\Lambda(kr)} \cm{K} \psi_{\Lambda(jr)} - \psi_{\Lambda(kr)} \cm{K} \psi_{\Lambda((k-1)r)} \| \le f(r), \qquad 1 \le j < k,
$$
and
$$
\| \psi_{\Lambda(kr)} \cm{K} - \psi_{\Lambda(kr)} \cm{K} \psi_{\Lambda((k-1)r)} \| \le f(r).
$$
From here it follows that
$$
\| \psi_k (\psi_\Lambda \cm{K})^k \psi_\Lambda - \psi_k \cm{K} \psi_{k-1} ... \cm{K} \psi_1 \cm{K} \psi_\Lambda \| \le (k-1) f(r) \|\cm{K}\|^{k-1}
$$
and
$$
\| \psi_k \cm{K}^k - \psi_k \cm{K} \psi_{k-1} ... \cm{K} \psi_1 \cm{K} \psi_\Lambda \| \le k f(r) \|\cm{K}\|^{k-1}.
$$
Therefore,
$$
\| \psi_{\Lambda(kr)} \big( \cm{K} \psi_{\Lambda(kr)} \big)^k - \psi_{\Lambda(kr)} \cm{K}^k \| \le 2k f(r) \|\cm{K}\|^{k-1},
$$
and replacing $r$ by $nr/k$, we get
$$
\| \psi_{\Lambda(nr)} \big( \cm{K} \psi_{\Lambda(kr)} \big)^k - \psi_{\Lambda(nr)} \cm{K}^k \| \le 2k f(r) \|\cm{K}\|^{k-1}
$$
for all $k=1,2,\ldots, n$.

Tacking into account all the estimations above, we can write
$$
\begin{array}{l}
\left\| \psi_{\Lambda(nr)} \big( 1-\psi_\Lambda \cm{K} \big)^{-1} \psi_\Lambda - \psi_{\Lambda(nr)} \big( 1- \cm{K} \big)^{-1} \right\| \le \\
\\
\le
\left\| \psi_{\Lambda(nr)} \big( \sum \limits_{k=n+1}^\infty (\psi_\Lambda \cm{K} )^k \big) \psi_\Lambda \right\| + \left\| \psi_{\Lambda(nr)} \big( \sum \limits_{k=n+1}^\infty \cm{K}^k \big) \right\| + \left\| \psi_{\Lambda(nr)} \sum \limits_{k=1}^n \big( (\psi_\Lambda \cm{K} )^k \psi_\Lambda - \cm{K}^k \big) \right\| \le \\
\\
\le \dfrac{2 \| \cm{K} \|^{n+1}}{1-\| \cm{K} \|} + 2 f(r) \sum \limits_{k=1}^n k \|\cm{K}\|^{k-1} < \dfrac{2 \| \cm{K} \|^{n+1}}{1-\| \cm{K} \|} + \dfrac{2 f(r)}{\big( 1 - \| \cm{K} \| \big)^2}.
\end{array}
$$
The right-hand side of the obtained inequality tends to zero as $n$ and $r$ increase. Putting $d = nr$, we can find a positive decreasing function $\varepsilon$ such that $\lim \limits_{d \to \infty} \varepsilon(d) = 0$ and
$$
\left\| \psi_{\Lambda(d)} \big( 1-\psi_\Lambda \cm{K} \big)^{-1} \psi_\Lambda - \psi_{\Lambda(d)} \big( 1- \cm{K} \big)^{-1} \right\| < \frac{\varepsilon(d)}{\|\delta\|}.
$$
From here it follows that
$$
\| \psi_{\Lambda(d)} \rho_\Lambda - \psi_{\Lambda(d)} \rho \| < \varepsilon(d),
$$
and, therefore, for any $I \subset \Lambda$ and any $x \in X_*^I$, we have
$$
\left| \rho_\Lambda (x) - \rho(x) \right| \le \| \psi_I \rho_\Lambda - \psi_I \rho \| < \varepsilon(d)
$$
as soon as $d(I,\Lambda^c) > d$.
\end{proof}

\begin{remark}
Let $\Phi$ be a converging pair interaction potential. Then the set $\Delta_1^{\Phi} = \{\Delta_t^{\bar x}, \bar x \in X^{t^c}, t \in \mathbb{Z}^d\}$ of functions
$$
\Delta_t^{\bar x}(x,u) = \sum \limits_{s \in t^c} (\Phi_{ts}(u \bar x_s) - \Phi_{ts}(x \bar x_s)), \qquad x,u \in X^t,
$$
forms a one-point transition energy field corresponding to the potential $\Phi$. It is not difficult to see that the elements of $\Delta_1^{\Phi}$ satisfy~\eqref{Delta-Phi-2}.

In particular, if $\Phi$ is a vacuum potential with the norm
$$
\|\Phi\| = \sup \limits_{t \in \mathbb{Z}^d} \sup \limits_{x \in X_*^t} \sum \limits_{s \in t^c} \sup \limits_{y \in X_*^s} \vert \Phi_{ts}(xy) \vert,
$$
then
$$
\|\Delta_1\| \le 2 \|\Phi\|, \qquad D \le \|\Phi\|,
$$
and Theorem~\ref{th-main} holds true if
$$
\frac{e^{ 2\|\Phi\|} N_X}{1 + e^{- 2\|\Phi\|} N_X} \left( 1 + 2 \big( 1 + 2 e^{2 \|\Phi\|} N_X \big) \big(\exp \{e^{\|\Phi\|} -1\} -1\big) \right) < 1.
$$
\end{remark}

\begin{remark}
It was shown in~\cite{DN19} (see also~\cite{KhN23}) that any probability distribution $P_\Lambda$  on $X^\Lambda$, $\Lambda \in W$, necessary has a Gibbs form:
$$
P_\Lambda(x) = \frac{e^{\Delta_\Lambda(x,u)}}{\sum \limits_{z \in X^\Lambda} e^{\Delta_\Lambda(z,u)}}, \qquad x \in X^\Lambda,
$$
where $u \in X^\Lambda$ and $\Delta_\Lambda$ is a transition energy function, that is, a function on $X^\Lambda \times X^\Lambda$ satisfying
$$
\Delta_\Lambda(x,u) = \Delta_\Lambda(x,z) + \Delta_\Lambda(z,u), \qquad x,u,z \in X^\Lambda.
$$
Correlation function $\rho_\Lambda$ for $P_\Lambda$ is defined as follows: for any $I \subset \Lambda$ and $x \in X_*^I$,
$$
\rho_\Lambda(x) = \sum \limits_{y \in X^{\Lambda \backslash I}} P_\Lambda(xy),
$$
and we come to the definition of finite-volume correlation functions used in this paper.
\end{remark}

\noindent \textbf{Funding.} The work was supported by the Science Committee of the Republic of Armenia in the frames of the research project 21AG-1A045.\\

\end{document}